\documentclass[12pt]{article}
\usepackage{amsmath,amsthm,amsfonts}
\usepackage{color}
\numberwithin{equation}{section} \allowdisplaybreaks

\newtheorem{theorem}{\sc Theorem}[section]
\newtheorem{lemma}[theorem]{\sc Lemma}
\newtheorem{proposition}[theorem]{\sc Proposition}
\newtheorem{corollary}[theorem]{\sc Corollary}
\newtheorem{definition}[theorem]{\sc Definition}
\newtheorem{example}[theorem]{\sc Example}
\newtheorem{remark}[theorem]{\sc Remark}

\newcommand{\bet}{\begin{theorem}}
\newcommand{\eet}{\end{theorem}}
\newcommand{\blm}{\begin{lemma}}
\newcommand{\elm}{\end{lemma}}
\newcommand{\bprop}{\begin{proposition}}
\newcommand{\eprop}{\end{proposition}}
\newcommand{\bcor}{\begin{corollary}}
\newcommand{\ecor}{\end{corollary}}
\newcommand{\bdf}{\begin{definition}\rm}
\newcommand{\edf}{\end{definition}}
\newcommand{\bp}{\begin{proof}}
\newcommand{\ep}{\end{proof}}
\newcommand{\bex}{\begin{example}\rm}
\newcommand{\eex}{\end{example}}
\newcommand{\bremark}{\begin{remark}\rm}
\newcommand{\eremark}{\end{remark}}

\begin{document}

\title {A trace formula for the index of B-Fredholm operators}

\author{ M. Berkani}

\date{}

 \maketitle
\vspace{-12mm}

\begin{abstract}

In this paper we define B-Fredholm elements in a Banach algebra $A$
modulo an ideal $J$ of $A.$ When a trace function is given on the
ideal $J,$ it generate an index for B-Fredholm elements. In the
case of a B-Fredholm operator $T$ acting on a Banach space, we
prove that its usual index $ind(T)$ is equal to the trace of the
commutator $ [T, T_0],$ where $T_0$ is a Drazin inverse of \,$T$
modulo the ideal of finite rank operators, extending a Fedosov's
trace formula for Fredholm operators \cite{BS}. In the case of a
semi-simple Banach algebra, we prove a punctured neighborhood
theorem for the index.
\end{abstract}

\renewcommand{\thefootnote}{}
\footnotetext{\hspace{-7pt} {\em 2010 Mathematics Subject
Classification\/}: Primary 47A53, 46H05
\baselineskip=18pt\newline\indent {\em Key words and phrases\/}:
B-Fredholm, Banach algebra, index,trace}

\section{Introduction}

Let $X$ be a Banach space and let $L(X)$ be the Banach algebra of
bounded linear operators acting on $X.$ In \cite{berkani-7}, we
have  introduced the class of linear bounded B-Fredholm operators.
If  $F_0(X)$ is the ideal of finite rank operators in $ L(X)$ and
 $\pi: L(X) \longrightarrow  A $ is the canonical projection,
from $L(X)$ onto the quotient algebra  $ A= L(X)/F_0(X),$  it is
well known by the Atkinson's theorem \cite [Theorem 0.2.2,
p.4]{BMW}, that $T \in L(X) $ is a Fredholm operator if and only
if its projection
 $\pi(T)$ in the   algebra $ A $ is invertible.  Similarly, in the
 following result, we have established an Atkinson-type theorem for B-Fredholm
operators.

\begin {theorem}\label{thm1}\cite [Theorem 3.4]{B-10}: Let $T \in L(X).$ Then $T$ is a B-Fredholm operator
if and only if $\pi(T)$ is Drazin invertible in the algebra $
L(X)/F_0(X).$
\end{theorem}

\noindent Motivated by this result,   we define in this paper
B-Fredholm elements in a semi-prime unital  Banach  algebra $A$
modulo an ideal $J$ of $A.$ Recall that  a Banach algebra $A$ is
called semi-prime if for $ u\in A,uxu= 0,$ for all $ x \in A$
implies that $u=0.$

\noindent An element $a \neq 0$ in a semiprime Banach algebra A is
called of rank one if there exists a linear functional $f_a$ on
$A$ such that $axa = f_a(x)a$ for all $x \in A.$

\bdf \label{def1}Let $A$ be  unital semi-prime Banach algebra and
let $J$ be an ideal of $A,$  and  $\pi: A\rightarrow A/J $  the
canonical projection. An element $a \in A$ is called  a B-Fredholm
element of  $A$ modulo the ideal $J$ if  its image $\pi(a)$ is
Drazin invertible in the quotient algebra $A/J.$ \edf

\noindent In a recent work \cite{CBS}, the authors, gave  in
\cite[Definiton 2.3]{CBS} a definition of   B-Fredholm elements in
Banach algebras. However, their definition does not englobe the
class of B-Fredholm operators, since the algebra $L(X)/F_0(X)$ is
not a Banach algebra. That's why in our definition, we consider
general algebras, not necessarily being Banach algebras, so it
includes also the case of the algebra $L(X)/F_0(X).$ While an
extensive study of B-Fredholm elements in a Banach algebras modulo
an ideal $J$ is being done in \cite{BMB},   we focus our attention
here on the properties of the index of such elements.

\noindent Recently in \cite{GR}, the authors studied Fredholm
elements in a semi-prime Banach algebra  modulo an ideal. When a
trace function is given on the ideal considered, they define also
the index of a Fredholm element using a trace formula
\cite[Definition 3.3]{GR}. This definition coincides with the
usual definition of the index for a Fredholm operator acting on a
Banach space.

\noindent In the second section of this paper, following the same
approach as in \cite{GR},  when a trace function is defined on the
ideal $J$ considered, we will define an index for  B-Fredholm
elements of the Banach algebra  $A$ modulo the ideal $J.$  Then in
the case of B-Fredholm operator $T$ acting on a Banach space $X,$
we prove our main result announced in the abstract, that is the
usual index $ind(T)$ is equal to $ \tau ([T, T_0])$ where $T_0$ is
a Drazin inverse of $T$ modulo the ideal of finite rank operators.

\noindent In the third section, we give a punctured neighborhood
theorem for the index B-Fredholm elements in a semi-simple Banach
algebra.  We will also establish a logarithmic rule of the index,
for two commuting prime B-Fredholm elements.

\section{Trace and Index}

In this section, $A$ will be a semi-prime, complex and unital
Banach algebra. We consider here B-Fredholm elements in $A$ modulo
an ideal $J,$ on which a trace function is given. Then we extend
the definition of the index \cite[Definition 3.3]{GR} to the class
of B-Fredholm elements of $A$ modulo $J.$ We will  show also that
the new definition of the index coincides with the usual one in
the case of B-Fredholm operator acting on a Banach space given in
\cite[Definition 2.3]{berkani-7}.

\bdf \label{def2}\cite[2.1., p.283] {GR} Let $J$ be an ideal in a
Banach algebra $A.$ A function $ \tau : J \rightarrow \mathbb{C},$
is called a trace on $J$ if :

1)$ \tau(p)= 1$ if $p \in J$ is an idempotent, that is $p^2=p,$
and $p$ of rank one,

 2) $\tau (a+b)= \tau (a) + \tau (b), $ for all $ a, b \in J,$

3) $\tau(\alpha a)= \alpha \tau(a),$ forall $ \alpha \in
\mathbb{C}$ and $ a \in J,$

4) $\tau(ab)= \tau(ba),$ for all $ a \in J$ and $ b \in  A.$

\edf

\bdf \label{def3}Let $\tau$ be a trace on an ideal $J$ of a Banach
algebra $A.$  The index of a B-Fredholm element $a \in A$ is
defined by: $$ \mathbf{i}(a)= \tau( aa_0 - a_0a)= \tau ([a,
a_0]),$$

\noindent where $ a_0 $ is a Drazin inverse of $a$ modulo the
ideal $I.$ \edf

\bet \label{thm2}The index of a B-Fredholm element $a \in A$ is
well defined and is independant of the   Drazin inverse $a_0 $  of
$a$ modulo the ideal $J.$

\eet

\bp Let $ a \in A$ be a B-Fredholm element in $A$ modulo $J.$ Then
$ \pi(a)$ is Drazin invertible in $A/J.$ If $\pi(a_0)$ is the
Drazin inverse of $\pi(a),$ then $\pi(a)\pi(a_0)= \pi(a_0)\pi(a).$
Hence $aa_0-a_0a \in J.$ If $\pi(a_0)= \pi(a_0'),$ then  $ a_0'-
a_0 \in J$ and $aa_0'-a_0'a- (aa_0 - a_0a)= a( a_0'-a_0)- (
a_0'-a_0)a.$ Since $ (a_0'-a_0) \in J,$  using the property 4) of
the trace, we obtain $ \tau(a( a_0'-a_0))= \tau (( a_0'-a_0)a).$
So the index is independent of the choice of the representative
$a_0.$ \ep

\bremark Clearly the definition  of the index of B-Fredholm
elements, extends \cite[Definition 3.3]{GR}, given for Fredholm
elements in Banach algebras, because the inverse modulo $J$ of a
Fredholm element $ a \in A$ is also its Drazin inverse modulo $J.$
\eremark

\noindent Now we prove in the next theorem, that the index defined
here for B-Fredholm elements in Banach algebras is equal to the
usual index of B-Fredholm operators acting on a Banach space. Here
the ideal $I$ is the ideal of finite rank operators acting on $X,$
and the trace $\tau$ is the  trace that is  defined for finite
rank operators as in \cite[Example 2.1]{GR}. If $T \in I$ and
$\{x_1, ..., x_n\}$ is a basis of its image $R(T),$ then $ T=
\Sigma_{i=1}^n x_i^\prime \bigotimes x_i,$ where  $x_1', ...,
x_n'$ are continuous linear functionals on $X$ such that $
T(x)=\Sigma_{i=1}^n x_i'(x)x_i ,$ for all $ x \in X.$ The trace of
$T$ is then defined by: $$ \tau(T)= \Sigma_{i=1}^n x_i'(x_i).$$

 \bet \label{thm3}Let $X$ be  a Banach space and let $T \in B(X)$
be a B-Fredholm operator. Then its usual index \,  $ind(T)$ is
equal to the trace $ \tau ([T, S])$ where $S$ is a Drazin inverse
of $T$ modulo the ideal of finite rank operators. \eet

\bp  Let $T$ be a B-Fredholm operator, then from \cite[Theorem
2.7] {berkani-7}  there exist two closed subspaces $M$ and $N$ of
$X $ such that
  $X = M\oplus N$  and:\\
\noindent i) $T(N)\subset N$ and  $T_{\mid N}$ is a  nilpotent operator,\\
\noindent ii) $T(M) \subset M$  and  $ T_{\mid M }$ is a Fredholm  operator. \\
Let $p \in \mathbb{N},$ such that $(T_{\mid N})^p=0.$ Then $ T^p=
(T_{\mid M})^p \oplus 0.$ From \cite[Theorem 2.7] {berkani-7},
$T^p$ is a B-Fredholm operator and  from \cite[Proposition 2.1]
{berkani-7} $ind( T^p)= p.ind(T)= p \, ind(T_{\mid M}). $ Let $
T_0$ be an inverse of $T_{\mid M}$ modulo the finite rank
operators on $M$ and $ S= T_0 \oplus 0$. Then $ S= T_0 \oplus 0$
is a Drazin inverse of $T$ and $ S^p= (T_0)^p \oplus 0$ is a
Drazin inverse of $T^p$ modulo the finite rank operators on $X.$
Moreover $T^p S^p - S^pT^p= [(T_{\mid M})^p T_0^p- T_0^p(T_{\mid
M})^p] \oplus 0. $

\noindent We observe  that $(T_{\mid M})^p T_0^p- T_0^p(T_{\mid
M})^p$ is of finite rank. So $T^p S^p - S^p T^p$ is also of finite
rank. Since $T_0^p$ is an inverse of the Fredholm  operator
$(T_{\mid M})^p$ modulo the finite rank on $T_{\mid M},$ then from
\cite [Example 3.2]{GR},we have $ \tau(T^p S^p- S^pT^p)= \tau(
(T_{\mid M})^p T_0^p- T_0^p(T_{\mid M})^p)= ind( (T_{\mid M})^p)=
p\, ind(T_{\mid M}))= p\,\,ind(T).$ Here $\tau,$ is the function
trace defined on finite rank operators acting on a Banach space as
in \cite[Examlpe 2.1] {GR}.

\noindent On another side, by \cite[Proposition 3.5]{GR} we have $
\tau((T_{\mid M})^p T_0^p- T_0^p(T_{\mid M})^p)= \mathbf{i}(
(T_{\mid M})^p)= p\,\mathbf{i}(T_{\mid M}) $ and so $ ind(T)=
\mathbf{i}(T_{\mid M}).$ Since, by definition of the index
$\mathbf{i},$ we have  $\mathbf{i}(T)= \tau(T S - S T)=
\tau((T_{\mid M} T_0 - T_0 T_{\mid M})\oplus 0) = \tau(T_{\mid M}
T_0 - T_0 T_{\mid M}) = \mathbf{i}(T_{\mid M}),$ we obtain $
ind(T)= \mathbf{i}(T).$

\ep

\bremark The trace formula for the index of B-fredholm operators
given in Theorem \ref{thm3}, extends Fedosov's trace formula for
the index of Fredholm operators, see\cite[p.10]{BS}

\eremark

\section{ Properties of the Index}

In this section, we will assume that $ A$ is a  semi-simple complex
unital Banach algebra, with unit $e,$ and the ideal $J$ is equal
to its socle. Recall that it is well known, that a semi-simple
Banach algebra is semi-prime. Then , it follows from \cite[BA2.4,
p. 103]{BMW} that an element $a$ in $A$ is invertible modulo $J$
if and only if it is invertible modulo the closure $\overline{J}$
of $J.$ In this case, if $p$ is any minimal idempotent in $A,$
that's a non zero idempotent such that $pAp= \mathbb{C}e,$ then
the operator $ \widehat{a}: Ap \rightarrow Ap,$ defined by
$\widehat{a}(x)= ax,$ is a Fredholm operator.

\noindent The element $a \in A $ is said to be of finite rank if
the operator $ \widehat{a}$  is an operator of finite rank. We
know from \cite[Theorem F.2.4] {BMW}, that the socle of $A$ is $
soc(A)= \{ x \in A\mid \widehat{x}$ {\text is of finite rank}\}.
Moreover, from \cite [Section 3]{AM}, a trace function is defined
on the socle by: $ \tau(a)= \Sigma_{\lambda \in \sigma(a)}
m(\lambda, a)\lambda,$ for an element $a$ of the socle of $A,$
where $ \sigma(a)$ is the spectrum of $a,$ and $m(\lambda, a)$ is
the algebraic multiplicity of $\lambda$ for $a.$

\noindent For more details about these notions from Fredholm
theory in Banach algebras, we refer the reader to \cite{BMW}.

\noindent In the following theorem, we will consider stability of B-Fredholmness under
small perturbations. Contrarily to the case of usual Fredholmness, we cannot except to preserve
B-Fredholmness under perturbation by  small norm elements. This can be easily in the case of  the algebra $L(X).$
Since  $0$ is a B-Fredholm operator on the Banach space $X,$ and there exists operators which are not
 B-Fredholm  (see \,\cite[Remark B]{B-12}), then we cannot have stability of
B-Fredholmness under perturbation by  small norm elements.

\bet \label{thm4}Let  $ a$ be a B-Fredholm element in $A$ modulo
$J.$ If $\lambda \in \mathbb{C}, \lambda \neq 0,$ and $|\lambda|$
is small enough, then $a-\lambda e $ is a Fredholm element of $A$
modulo $J$ and $ \mathbf{i}(a -\lambda e)= \mathbf{i}( a).$ \eet

\bp Assume that $a$ is a B-Fredholm element in $A$ modulo $J.$
Then $ \pi(a)$ is Drazin invertible in $A/J.$ Let $\Pi:
A\rightarrow A/\overline{J}, $  be  the canonical projection. Then
$ \Pi(a)$  is Drazin invertible in $ A/ \overline{J}.$ As $ A/
\overline{J}$ is a Banach algebra, then $0$ is isolated in the
spectrum of $\Pi(a)$ in the Banach algebra $ A/ \overline{J}.$
Thus if $|\lambda|$ is small enough and $\lambda \neq 0,$ then
$\Pi(a-\lambda e) $ is invertible in $ A/ \overline{J}.$ From our
hypothesis on the ideal $J,$ it follows that $\pi(a-\lambda e) $
is invertible in $ A/J.$ So $a-\lambda e $ is a Fredholm element
in $A$ modulo $J.$

\noindent Moreover,we have from \cite[Theorem F.2.6]{BMW}, if $ p$
is a minimal idempotent in $A,$ then the operator $
\widehat{a-\lambda I}: Ap \rightarrow  Ap,$ defined by
$\widehat{a-\lambda I}(x)= (a-\lambda I)x,$  is a Fredholm
operator on the Banach space $Ap.$ From \cite[Theorem 3.17]{GR},
it follows that $ ind( \widehat{a- \lambda I})= \mathbf{i}( a-
\lambda I).$

\noindent On another side, from \cite [ Remark, iii)]{berkani-7},
we have $ ind( \widehat{a})= ind( \widehat{a- \lambda I}), $ for $
|\lambda| ,$ small enough. Hence $ind( \widehat{a})= \mathbf{i}(
a- \lambda I),$ for $ |\lambda| $ small enough.

\blm \label{lem1} If $ p$ is a minimal idempotent in $A,$ then the
operator $ \widehat{a}: Ap \rightarrow  Ap,$ defined by
$\widehat{a}(x)= ax,$ is a B-Fredholm operator and
$ind(\widehat{a})=\mathbf{i}(a).$
 \elm

\bp Since $ a$ is a B-Fredholm element in $A$ modulo $J,$ then $a$
is Drazin invertible in $A$ modulo $J.$   From \cite[Theorem
F.2.4]{BMW}, we know that  J is exactly the set of elements $x$ of
$A$ such that $ \widehat{x}$ is an operator of finite rank. Then
$\widehat{a}$ is a Drazin invertible  operator modulo the ideal of
finite rank on $Ap.$ Thus from Theorem \ref{thm1}, $\widehat{a}:
Ap \rightarrow Ap $ is a B-Fredholm operator. Let $b \in A$ be a
Drazin inverse of $a$ modulo $J,$ then $\widehat{b}$ is a Drazin
inverse of $\widehat{a}$ modulo the ideal of finite rank on the
Banach space $Ap.$ From Theorem \ref{thm3}, we have
$\mathbf{i}(\widehat{a})=\tau(\widehat{a}\widehat{b}-\widehat{b}\widehat{a})=
\tau(\widehat{ab-ba})=ind(\widehat{a}).$ Here $\tau$ stands for
the trace of the finite rank operators on the Banach space $Ap.$
\ep

\noindent Using Lemma \ref{lem1}, and since $ab-ba \in J,$ it is of finite trace and
$\tau(ab-ba)= \tau(\widehat{ab-ba}).$ Thus $ \mathbf{i}(a)=
\mathbf{i}( a- \lambda I).$ \ep

As seen in \cite{B-16}, the product of  two
B-Fredholm operators in $L(X),$  even it is  a B-Fredholm operator, does not  have in general  its
 index equal to the sum of the indexes of the operators involved in the product,
 unless the two operators  satisfies a commuting  Bezout identity, as proved in \cite[Theorem 1.1]{B-16}. Here we obtain a similar result for the product of B-Fredholm elements.

\bprop Let  $ a_1, a_2$ be B-Fredholm elements in $A$ modulo $J,$
and let $ \lambda \in  \mathbb{C}.$ \\
 \noindent  i) If $ a_1, a_2, u_1, u_2 \in A$  are two by two commuting elements in $A$  such that
 $ u_1a_1+ u_2a_2=e,$ then $a_1a_2$ is a B-Fredholm
element in $A$ modulo $J,$ and $ \mathbf{i}( a_1a_2)=
\mathbf{i}(a_1)+  \mathbf{i}(a_2).$
 In particular, if $ \lambda \neq 0,$ then $\lambda a_1$ is a B-Fredholm element of $A$ modulo $J$ and  $ \mathbf{i}(
\lambda a_1)= \mathbf{i}( a_1).$\\
\noindent ii) If $j$ is an element of $J,$ then $a_1+j$ is a
B-Fredholm element in $A$ modulo $J$ and $ \mathbf{i}( a_1 +j)=
\mathbf{i}( a_1).$ \eprop

 \bp i) From \cite[Proposition 2.6]{B-10},
 it follows that  $a_1a_2$ is a B-Fredholm
element in $A$ modulo $J.$ Moreover as $ \widehat{u_1}
\widehat{a_1}+ \widehat{u_2}\widehat{a_2}= I,$ and $
\widehat{a_1}, \widehat{a_2}, \widehat{u_1}, \widehat{u_2}, $ are
two by two commuting operators, then  from \cite [Theorem
1.1]{B-16}, we have $ ind (\widehat{a_1}\widehat{a_2})= ind(
\widehat{a_1}) + ind( \widehat{a_2}),$ and  Lemma \ref{lem1}
implies that  $ \mathbf{i} (\widehat{a_1}\widehat{a_2})=
\mathbf{i}( \widehat{a_1}) + \mathbf{i}( \widehat{a_1}). $ Taking
$ a_2= \lambda e,$ we obtain $ \mathbf{i}( \lambda a_1)=
\mathbf{i}( a_1)$ because $ \mathbf{i}( \widehat{\lambda e})= ind(
\lambda I)= 0.$

\noindent ii) If $j$ is an element of $J,$ then $ \pi(a_1+j)=
\pi(a_1).$ So  $ a_1+j$ is a B-Fredholm element in $A$ modulo $J.$
Moreover  if $|\lambda|$ is small enough and $\lambda \neq 0,$
then from Theorem \ref{thm4}, $a_1 - \lambda e + j$ is a Fredholm
element and $ \mathbf{i}( a_1 +j)= \mathbf{i}( a_1 - \lambda I +
j)).$ Using \cite[Proposition 3.7, i]{GR}, we obtain $\mathbf{i}(
a_1 - \lambda I + j))= \mathbf{i}( a_1 - \lambda I)= \mathbf{i}(
a_1).$
 \ep

 \baselineskip=12pt
\bigskip
\vspace{10 mm }
 \baselineskip=12pt
\bigskip

{\small
\noindent Mohammed Berkani,\\
 \noindent Department of Mathematics,\\
 \noindent Science faculty of Oujda,\\
\noindent University Mohammed I,\\
\noindent Laboratory LAGA, \\
\noindent Morocco\\
\noindent berkanimo@aim.com,\\

\end{document}